\documentclass[12pt,halfline,a4paper]{ouparticle}
\newtheorem{theorem}{Theorem}[section]

\begin{document}

\title{On the Zagreb indices of the line graph and line cut-vertex graph of subdivision of unicyclic graphs}
\author{%
\name{H. M. Nagesh}
\address{Department of Science and Humanities, \\
PES University, Bangalore, India}
\email{hmnagesh1982@gmail.com}}

\abstract{The first Zagreb index $M_{1}(G)$ is equal to the sum of squares of the degrees of the vertices, and the second Zagreb index $M_{2}(G)$  is equal to the sum of the products of the degrees of pairs of adjacent vertices of the underlying molecular graph $G$. This paper aims to investigate the Zagreb indices and coindices of the line graph and line cut-vertex graph of a class of unicyclic graphs called cycle-star graph using the notion of subdivision. 

\vskip1em \noindent \textbf{2010 AMS Classification:} 05C07, 05C35, 05C90.

} 
\date{}
\maketitle
\label{sec1}
\section{Introduction} \label{sec:Intr}
Let $G=(V,E)$ be a simple graph with vertex set $V(G)=\{v_1,v_2,\ldots,v_n\}$ and edge set $E(G)$, where $|V(G)|=n$ and $|E(G)|=m$. These two basic parameters $n$ and $m$ are called the $order$ and $size$ of $G$, respectively, and they are essential for characterizing and quantifying the basic properties of a graph. When referring to an edge that connects two vertices $u$ and $v$, it is common to use the notation $uv$ to represent the edge between $u$ and $v$. The degree of a vertex in a graph $G$ is a fundamental concept in graph theory. The $degree$ of a vertex $v$, written $d_{G}(v)$, is defined as the number of edges in the graph $G$ that are incident with vertex $v$.  

The first and second Zagreb indices are among the oldest and most studied topological indices in the field of chemical graph theory. These indices are used to describe the topological structure of chemical compounds, particularly organic molecules. These two indices first appeared in \cite{4}, and were elaborated in \cite{5}. The main properties of $M_1(G)$ and $M_2(G)$ were summarized in \cite{7,14}. 

Certainly, the first Zagreb index $M_1(G)$ and the second Zagreb index $M_2(G)$ of a graph $G$ are defined, respectively, as
\begin{equation*}
M_{1}=M_{1}(G)=\displaystyle \sum_{v \in V(G)} d_{G} (v)^2=\displaystyle \sum_{uv \in E(G)} [d_{G}(u)+d_{G}(v)]
\end{equation*}
\begin{equation*}
M_{2}=M_{2}(G)=\displaystyle \sum_{uv \in E(G)} d_{G}(u) \cdot d_{G}(v)
\end{equation*}

Indeed, during the past decades, numerous results and studies have been conducted concerning the Zagreb indices ($M_1$ and $M_2$) and their applications in various fields, particularly in chemical graph theory and chemo-informatics \cite{2,3,6,8,9}. Some historical details about the Zagreb indices can be seen in \cite{10}. 

In 2008, Do\v{s}li\'c introduced the concept of the first Zagreb coindex as an extension or complement to the existing Zagreb indices. The first Zagreb coindex, denoted by $\overline{M_1}$, is a topological index that is defined based on the same principles as the original Zagreb indices. 
\begin{equation} 
\overline{M_1}=\overline{M_{1}}(G)=\displaystyle \sum_{uv \notin E(G)} [d_{G}(u)+d_{G}(v)]
\end{equation}
The second Zagreb coindex, denoted by $\overline{M_2}$, is defined analogously to the second Zagreb index ($M_2$). 
\begin{equation} 
\overline{M_2}=\overline{M_{2}}(G)=\displaystyle \sum_{uv \notin E(G)} d_{G}(u) \cdot d_{G}(v)
\end{equation}
In expressions (1) and (2) for the first and second Zagreb coindex, it is indeed assumed that the vertices $u$ and $v$ represent distinct vertices in the graph $G$, meaning that $u \neq v$. 

The \emph{subdivision graph} of a graph $G$, denoted by $S(G)$, is a graph that is constructed from $G$ by introducing additional vertices along the edges of $G$. This process involves replacing each edge in G with a path of length 2, effectively subdividing each edge. 

Jelena Sedlar \cite{20} introduced the notion of the cycle-star graph. A \emph{cycle-star graph}, denoted by $CS_{k,n-k}$, is a type of graph that consists of two components: a cycle of length $k$ and $n-k$ leaf vertices connected to the same vertex of the cycle. 

So, all cycle-star graphs $CS_{k,n-k}$ are unicyclic graphs (i.e., connected graphs containing exactly one cycle). An example of a cycle-star graph is shown in Figure 1.
\newpage
\vspace{5mm}
\begin{center}
\includegraphics[width=9cm]{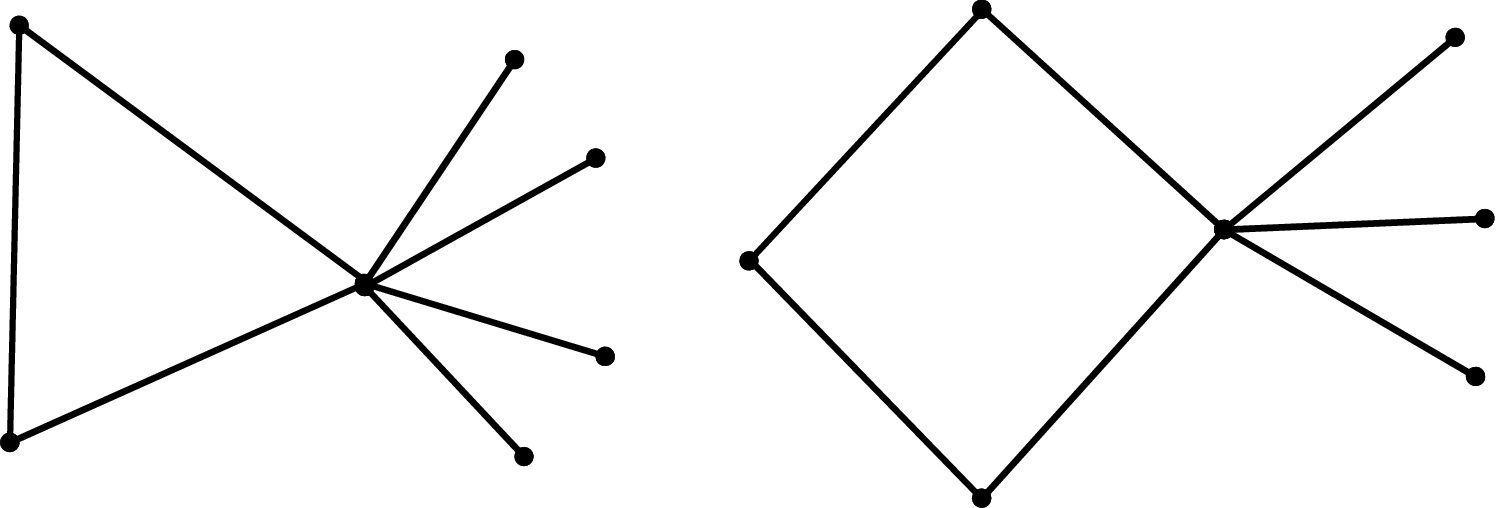}
\end{center}
\vspace{5mm}
\begin{center}
Figure 1. The cycle-star graphs $CS_{3,4}$ and $CS_{4,3}$ 
\end{center}

Unicyclic graphs are indeed of significant interest in the field of chemistry because they can represent a wide range of chemical structures, especially in the context of organic chemistry. Unicyclic graphs are used to model and analyze the structures of organic compounds, which are compounds primarily composed of carbon atoms bonded to other elements like hydrogen, oxygen, nitrogen, and so on. 

Recently, Zagreb indices of unicyclic graphs attracted much attention. Studies along this line include general multiplicative Zagreb indices of unicyclic graphs \cite{13}, Zagreb eccentricity indices of unicyclic graphs \cite{21}, Maximal hyper-Zagreb index of unicyclic graphs with a given order and matching number \cite{22}. However, the studies on the Zagreb indices of the intersection graph on the vertex set of unicyclic graphs were not attempted. 

In this paper, we study the line graph and line cut-vertex graph of the subdivision graph of a class of unicyclic graphs called the cycle-star graph; and calculate the Zagreb indices and coindices of these graphs. To find more details and in-depth information about Zagreb indices and coindices of graphs using graph operators, readers are referred to the articles \cite{15,16,17}. Notations and definitions not introduced here can be found in \cite{11}.

The present paper is organized as follows. In Section 2, we find the Zagreb indices and coindices of the line graph of the subdivision graph of the cycle-star graph $CS_{k,n-k}$. In Section 3, we find the Zagreb indices and coindices of the line cut-vertex graph of the subdivision graph of the cycle-star graph $CS_{k,n-k}$. Finally, Section 4 gives the conclusion and directions for further research. The paper is completed with a list of references.

\section{Zagreb indices of the line graph of the subdivision graph of the cycle-star graph $CS_{k,n-k}$} 
In this section, the Zagreb indices and coindices of the line graph of the subdivision graph of the cycle-star graph 
are calculated.

In particular, unless otherwise specified, the parameters $k \geq 3$ and $n-k \geq 1$ are considered for every cycle-star graph $CS_{k,n-k}$.

There are many graph operators (or graph valued functions) with which one can construct a new graph from a given graph, such as line graphs, line cut-vertex graphs; total graphs; middle graphs; and their generalizations. 

The \emph{line graph} of a graph $G$, written $L(G)$, is the graph whose vertices are the edges of $G$, with two vertices of $L(G)$ adjacent whenever the corresponding edges of $G$ have a vertex in common. In \cite{19}, the Zagreb indices and coindices of the line graphs of the subdivision graphs were studied.

In the next theorem, the line graph of the subdivision graph of the cycle-star graph is determined.

\begin{theorem}
Let $G$ be the line graph of the subdivision graph of the cycle-star graph $CS_{k,n-k}$. Then $M_1(G)=n^3+(6-3k)n^2+(3k^2-12k+13)n-k^3+6k^2-5k$ 
and \begin{eqnarray*}
M_2(G)=\frac{1}{2} \left( n^4+(7-4k)n^3+(6k^2-21k+20)n^2-(4k^3-21k^2+40k-32)n\right)+ \\
\frac{1}{2} \left( k^4-7k^3+20k^2-16k \right).
\end{eqnarray*}
\end{theorem}
Proof. Let $G$ be the line graph of the subdivision graph of the cycle-star graph $CS_{k,n-k}$, that is, $G=L(S(CS_{k,n-k}))$. The subdivision graph $S(CS_{k,n-k})$ contains $2n$ vertices and $2n$ edges so that the line graph of $S(CS_{k,n-k})$ contains $2n$ vertices, out of which $2k-2$ vertices of are of degree $2$; $n-k+2$ vertices are degree $n-k+2$, and the remaining $n-k$ vertices are of degree $1$. Thus,
 \begin{align*}
M_1(G) & = 4(2k-2)+(n-k+2)(n-k+2)^2+(n-k) \\
& = n^3+(6-3k)n^2+(3k^2-12k+13)n-k^3+6k^2-5k.
\end{align*} 

To find $M_2(G)$, we first find the size of $G$. The size of $G$ is
 \begin{align*}
|E(G)| & = \frac{(n-k+2)(n-k+1)}{2}+2k-1+n-k 
\end{align*} 
But,
\begin{equation*}
 \frac{(n-k+2)(n-k+1)}{2}=\frac{1}{2}[n^2+k^2+3n-3k-2nk+2]
\end{equation*}
Hence,
\begin{align*}
|E(G)|  & = \frac{1}{2}[n^2+k^2+5n-2nk-k]
\end{align*} 
In other words, $E(G)$ contains $2k-3$ edges whose end vertices have degree $2$; $2$ edges whose end vertices have degree $2$ and $n-k+2$; $n-k$ edges whose end vertices have degree $1$ and $n-k+2$; and the remaining $\frac{(n-k+2)(n-k+1)}{2}$ edges whose end vertices have degree $n-k+2$.
Thus, 
\begin{equation*}
M_2(G)=(8k-12)+(4n-4k+8)+(n-k)(n-k+2)+\frac{(n-k+2)(n-k+1)}{2}(n-k+2)^2
\end{equation*} 
But, 
\begin{align*}
(n-k)(n-k+2)=n^2+k^2-2nk+2n-2k; (n-k+2)^2 & = n^2+k^2-2nk+4n-4k+4 
\end{align*} 
\begin{eqnarray*}
\frac{(n-k+2)(n-k+1)}{2}(n-k+2)^2  = \frac{1}{2} \left(n^4+(7-4k)n^3+(6k^2-21k+18)n^2 \right)-\\
\frac{1}{2} \left((4k^3-21k^2+36k-20)n+k^4-7k^3+18k^2-20k+8 \right)
\end{eqnarray*}
Hence,
\begin{eqnarray*}
M_2(G)=\frac{1}{2} \left( n^4+(7-4k)n^3+(6k^2-21k+20)n^2-(4k^3-21k^2+40k-32)n\right)+ \\
\frac{1}{2} \left( k^4-7k^3+20k^2-16k \right)
\end{eqnarray*}
An example of Theorem 2.1 is shown in Figure 3. 
\vspace{3mm}
\begin{center}
\includegraphics[width=15cm]{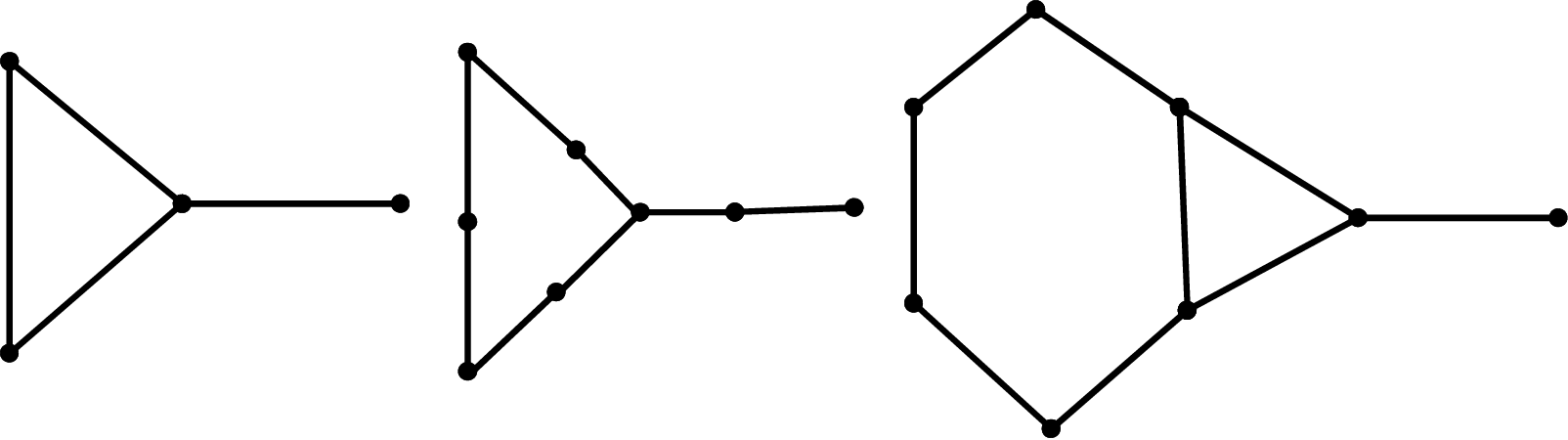}
\end{center}
\begin{center}
Figure 3. The cycle-star graph $CS_{3,1}$; $S(CS_{3,1})$, $L(S(CS_{3,1}))$
\end{center}

For $G=L(S(CS_{3,1}))$ of Figure 3, using Theorem 2.1, 
\begin{eqnarray*}
M_{1}(G)=64+16(6-9)+4(27-36+13)-27+54-15=44.
\end{eqnarray*}
\begin{eqnarray*}
M_2(G)=\frac{1}{2} \left(256+64(7-12)+16(54-63+20)-4(108-189+120-32) \right)+\\
 \frac{1}{2} \left(81-189+180-48 \right) \\
={} 54.
\end{eqnarray*}
We now give the expressions for the first and second Zagreb coindices of the line graph of the subdivision graph of the cycle-star graph.

Gutman et al. in \cite{8} established a complete set of relations between the first and second Zagreb index and coindex of a graph as follows: 

\begin{theorem}
Let $G$ be a graph with $n$ vertices and $m$ edges. Then  \begin{center}
$\overline{M_{1}}(G)=2m(n-1)-M_1(G)$  \end{center}
\end{theorem}          

\begin{theorem}
Let $G$ be a graph with $n$ vertices and $m$ edges. Then \begin{center}
$\overline{M_{2}}(G)=2m^2-\frac{1}{2} M_1(G)-M_2(G)$ \end{center}
\end{theorem}
Using Theorem 2.2 and Theorem 2.3, we have the following result.
\begin{theorem}
Let $G$ be the line graph of the subdivision graph of the cycle-star graph $CS_{k,n-k}$. Then $\overline{M_1}(G)=n^3+(3-k)n^2-(k^2-12k+18)n+k^3-7k^2+6k$.
\end{theorem}
Proof. Let $G$ be the line graph of the subdivision graph of the cycle-star graph $CS_{k,n-k}$. The order and size of $G$ are $2n$ and $\frac{1}{2}[n^2+k^2+5n-2nk-k]$, respectively. By Theorem 2.1, $M_1(G)=n^3+(6-3k)n^2+(3k^2-12k+13)n-k^3+6k^2-5k$. Then Theorem 2.2 implies that 
\begin{eqnarray*}
\overline{M_1}(G) = (n^2+k^2+5n-2nk-k)(2n-1)-(n^3+(6-3k)n^2+\\
(3k^2-12k+13)n-k^3+6k^2-5k)
\end{eqnarray*}
Since
\begin{equation*}
(n^2+k^2+5n-2nk-k)(2n-1)=2n^3+2nk^2+10n^2-4n^2k-2nk-n^2-k^2-5n+2nk+k,
\end{equation*}
\begin{align*}
\overline{M_1}(G) & = n^3+(3-k)n^2-(k^2-12k+18)n+k^3-7k^2+6k.
\end{align*}

\begin{theorem}
Let $G$ be the line graph of the subdivision graph of the cycle-star graph $CS_{k,n-k}$. Then 
$\overline{M_2}(G)=n^3+\frac{n^2}{2}(2k-1)-\frac{n}{2}(8k^3-32k^2+38k-19)+\frac{1}{2}(6k^3-25k^2+21k)$. 
\end{theorem}
Proof. By Theorem 2.1, $M_1(G)=n^3+(6-3k)n^2+(3k^2-12k+13)n-k^3+6k^2-5k$; \\
$M_2(G)=\frac{1}{2}[n^4+(7-4k)n^3+(6k^2-21k+20)n^2-(4k^3-21k^2+40k-32)n+k^4-7k^3+20k^2-16k]$. Then Theorem 2.3 implies that 
\begin{align*}
\overline{M_2}(G) & = \frac{1}{2} \left(n^2+k^2+5n-2nk-k \right)^2-\frac{1}{2}M_1(G)-M_2(G).
\end{align*}
Since,
\begin{equation*}
(n^2+k^2+5n-2nk-k)^2=n^4+(10-4k)n^3+(6k^2-22k+25)n^2-(4k^3-14k^2+10k)n+k^4-2k^3+k^2,
\end{equation*}
\begin{eqnarray*}
\overline{M_2}(G)=n^3+\frac{n^2}{2}(2k-1)-\frac{n}{2}(8k^3-32k^2+38k-19)+\frac{1}{2}(6k^3-25k^2+21k).
\end{eqnarray*}

\section{Zagreb indices of the line cut-vertex graph of the subdivision graph of the cycle-star graph $CS_{k,n-k}$}
In this section, the Zagreb indices and coindices of the line cut-vertex graph of the subdivision graph of the cycle-star graph are calculated. 

Kulli et al. \cite{12} introduced the notion of line cut-vertex graph of a graph. The \emph{line cut-vertex graph} of a graph $G$, written $L_{c}(G)$, is the graph whose vertices are the edges and cut-vertices of $G$, with two vertices of $L_{c}(G)$ adjacent whenever the corresponding edges of $G$ have a vertex in common; or one corresponds to an edge $e_i$ of $G$ and the other corresponds to a cut-vertex $c_j$ of $G$ such that $e_i$ is incident with $c_j$. 

Clearly, $L(G) \subseteq L_{c}(G)$, where $\subseteq$ is the subgraph notation. 
An example of a graph $G$ and its $L_{c}(G)$ is shown in Figure 2.
\vspace{5mm}
\begin{center}
\includegraphics[width=9cm]{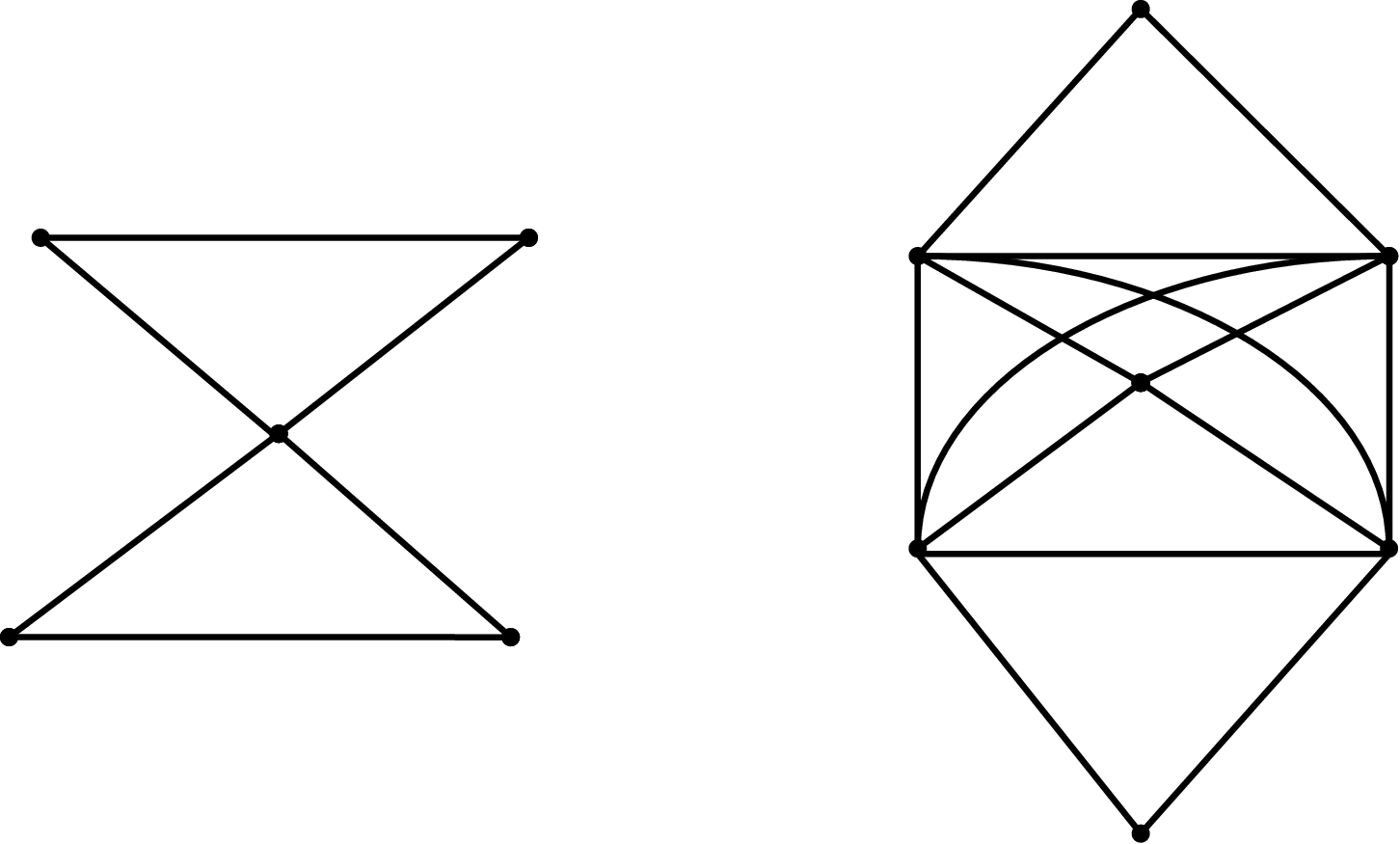}
\end{center}
\begin{center}
Figure 2. A graph $G$ and its line cut-vertex graph $L_{c}(G)$
\end{center} 
In the next theorem, the line cut-vertex of the subdivision graph of the cycle-star graph is determined.
\begin{theorem}
Let $G$ be the line cut-vertex graph of the subdivision graph of the cycle-star graph $CS_{k,n-k}$. Then,\\ $M_1(G)=n^3+(11-3k)n^2+(3k^2-22k+40)n-k^3+11k^2-32k+14$ 
and \begin{equation*}
\begin{aligned}
M_2(G)=\frac{1}{2}(n^4+k^4+(13-4k)n^3+(6k^2-37k+56)n^2-(4k^3-35k^2+96k-80)n)- \\
\frac{11}{2}k^3+22k^2-34k+9
\end{aligned}
\end{equation*}
\end{theorem}
Proof. Let $G$ be the line cut-vertex graph of the subdivision graph of the cycle-star graph $CS_{k,n-k}$, that is, $G=L_{c}(S(CS_{k,n-k}))$. The subdivision graph $S(CS_{k,n-k})$ contains $2n$ vertices and $2n$ edges so that the line cut-vertex graph of $S(CS_{k,n-k})$ contains $3n-k+1$ vertices, out of which $2k-2$ vertices of are of degree $2$; $2(n-k)$ vertices are degree $2$; $2$ vertices are of degree $n-k+3$; $1$ vertex is of degree $n-k+2$, and the remaining $n-k$ vertices are of degree $n-k+4$. Thus,
 \begin{align*}
M_1(G) & = 4(2k-4)+4(2n-2k)+2(n-k+3)^2+(n-k+2)^2+(n-k)(n-k+4)^2
\end{align*} 
But, 
\begin{align*}
2(n-k+3)^2 & = 2n^2+2k^2-4nk+12n-12k+18
\end{align*} 
\begin{align*}
(n-k+2)^2 & = n^2+k^2-2nk+4n-4k+4
\end{align*} 
\begin{align*}
(n-k+4)^2 & = n^2+k^2-2nk+8n-8k+16
\end{align*} 
Hence, 
\begin{eqnarray}\nonumber
M_1(G)=n^3+(11-3k)n^2+(3k^2-22k+40)n-k^3+11k^2-32k+14.
 \end{eqnarray}
Now, the size of $G$ is
 \begin{align*}
|E(G)| & = \frac{(n-k+3)(n-k+2)}{2}+3(n-k)+2k-1 \\
& = \frac{1}{2}(n^2+k^2+5n-5k-2nk+6)+3n-3k+3k-1 \\
& = \frac{1}{2}(n^2+k^2+11n-7k-2nk+4).
\end{align*} 
In other words, $E(G)$ contains $2k-3$ edges whose end vertices have degree $2$; $n-k$ edges whose end vertices have degree $2$;
$2$ edges whose end vertices have degree $2$ and $n-k+3$; $2(n-k)$ edges whose end vertices have degree $2$ and $n-k+4$; $2$ edges whose end vertices have degree $n$ and $n-1$; $n-k$ edges whose end vertices have degree $n-1$ and $n-k+4$; $1$ edge whose end vertices have degree $n-k+3$; $2(n-k)$ edges whose end vertices have degree $n-k+3$ and $n-k+4$; and the remaining $\frac{(n-k)(n-k-1)}{2}$ edges whose end vertices have degree 
$n-k+4$. Thus,  
\begin{equation*}
\begin{aligned}
M_2(G) = {} 4(2k-3)+4(n-k)+4(n-k+3)+2(n-k)(2n-2k+8)+2n(n-1)+\\
 (n-k)(n-1)(n-k+4)+(n-k+3)^2+2(n-k)(n-k+3)(n-k+4)+  \\
 \frac{(n-k)(n-k-1)}{2}(n-k+4)^2
\end{aligned}
\end{equation*}
But,
\begin{equation*}
2(n-k)(2n-2k+8)=4n^2+4k^2+16n-16k-8nk
\end{equation*}
\begin{equation*}
(n-k)(n-1)(n-k+4)=n^3+(3-2k)n^2+(k^2-2k-4)n-k^2+4k
\end{equation*}
\begin{equation*}
(n-k+3)^2=n^2+k^2+6n-6k-2nk+9
\end{equation*}
\begin{equation*}
2(n-k)(n-k+3)(n-k+4)=2n^3-2k^3+(14-6k)n^2+(6k^2-28k+24)n+14k^2-24k
\end{equation*}
\begin{equation*}
\frac{(n-k)(n-k-1)}{2}=\frac{1}{2}(n^2+k^2-n+k-2nk)
\end{equation*}
\begin{equation*}
\begin{aligned}
\frac{(n-k)(n-k-1)}{2}(n-k+4)^2=\frac{1}{2}(n^4+k^4+(7-4k)n^3+(6k^2-21k+8)n^2-\\
(4k^3-21k^2+16k+16)n-7k^3+8k^2+16k)
\end{aligned}
\end{equation*}
Hence,
\begin{equation*}
\begin{aligned}
M_2(G)=\frac{1}{2}(n^4+k^4+(13-4k)n^3+(6k^2-37k+56)n^2-(4k^3-35k^2+96k-80)n)- \\
\frac{11}{2}k^3+22k^2-34k+9.
\end{aligned}
\end{equation*}
\vspace{7mm}
An example of Theorem 3.1 is shown in Figure 4. 
\vspace{5mm}
\begin{center}
\includegraphics[width=15cm]{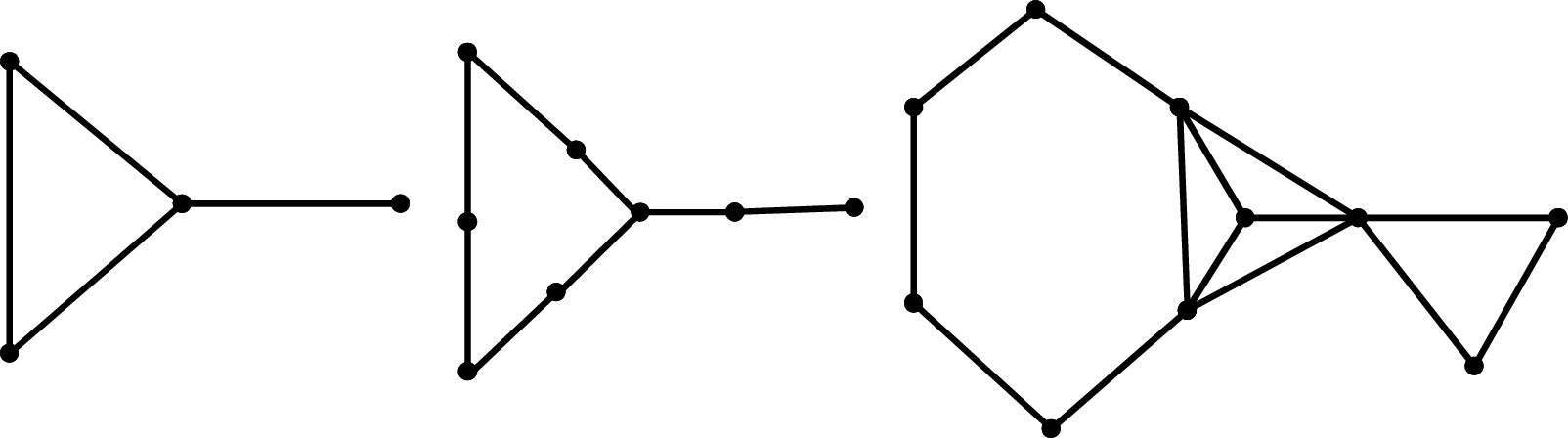}
\end{center}
\begin{center}
Figure 4: The cycle-star graph $CS_{3,1}$; $S(CS_{3,1})$, $L_{c}(S(CS_{3,1}))$
\end{center}
For $G=L_{c}(S(CS_{3,1}))$ of Figure 4, using Theorem 3.1, 
\begin{center}
$M_{1}(G)=64+16(11-9)+4(27-66+40)-27+99-96+14=90$.
\end{center}
\begin{align*}
\begin{aligned}
M_2(G)= {} \frac{1}{2}[256+64(13-12)+16(54-111+56)-4(108-315+288-80)]+\frac{81}{2}-\\
\frac{(11)(27)}{2}+198-102+9 \\
={} 147.
\end{aligned}
\end{align*}
We now give the expressions for the first and second Zagreb coindices of the line cut-vertex graph of the subdivision graph of the cycle-star graph using Theorem 2.2 and Theorem 2.3. \newpage
\begin{theorem}
Let $G$ be the line cut-vertex graph of the subdivision graph of the cycle-star graph $CS_{k,n-k}$. Then 
$ \overline{M_1}(G)=2n^3+(22-4k)n^2+(2k^2-10k-28)n-4k^2+28k-14$.
\end{theorem}
Proof. The order and size of $G$ are $3n-k+1$ and the size of $\frac{1}{2}(n^2+k^2+11n-7k-2nk+4)$, respectively. By Theorem 3.1, 
\begin{center}
$M_1(G)=n^3+(11-3k)n^2+(3k^2-22k+40)n-k^3+11k^2-32k+14$.
\end{center}
Then Theorem 2.2 implies that 
\begin{align*}
\begin{aligned}
\overline{M_1}(G) = {} (n^2+k^2+11n-7k-2nk+4)(3n-k)-((n^3+(11-3k)n^2+(3k^2-22k+40)n- \\
k^3+11k^2-32k+14))\\
= (3n^3+(33-7k)n^2+(5k^2-32k+12)n-k^3+7k^2-4k)-((n^3+(11-3k)n^2+\\
(3k^2-22k+40)n-k^3+11k^2-32k+14))\\
= 2n^3+(22-4k)n^2+(2k^2-10k-28)n-4k^2+28k-14.
\end{aligned}
\end{align*}

\begin{theorem}
Let $G$ be the line cut-vertex graph of the subdivision graph of the cycle-star graph $CS_{k,n-k}$. Then \begin{eqnarray*}
 \overline{M_2}(G)= \frac{1}{2} \left( 8n^3+(62-18k)n^2+(12k^2-52k-32)n \right)-\frac{1}{4}k^4-\frac{15}{4}k^3+12k^2+5k-\frac{7}{2} \end{eqnarray*}
\end{theorem}
Proof. By Theorem 3.1, $M_1(G)=n^3+(11-3k)n^2+(3k^2-22k+40)n-k^3+11k^2-32k+14$ 
and \begin{equation*}
\begin{aligned}
M_2(G)=\frac{1}{2}(n^4+k^4+(13-4k)n^3+(6k^2-37k+56)n^2-(4k^3-35k^2+96k-80)n)- \\
\frac{11}{2}k^3+22k^2-34k+9
\end{aligned}
\end{equation*}
Then Theorem 2.3 implies that 
\begin{align*}
\overline{M_2}(G) & = \frac{1}{2} \left(n^2+k^2+11n-7k-2nk+4 \right)^2-\frac{1}{2}M_1(G)-M_2(G).
\end{align*}
Since,
\begin{eqnarray*}
\left(n^2+k^2+11n-7k-2nk+4 \right)^2=\frac{1}{2}\left(n^4+k^4+(22-4k)n^3+(6k^2-58k+129)n^2 \right) \\
-\frac{1}{2}\left((4k^3-50k^2+170k-88)n+14k^3-57k^2+56k-16\right),
\end{eqnarray*}
\begin{eqnarray*}
 \overline{M_2}(G)= \frac{1}{2} \left( 8n^3+(62-18k)n^2+(12k^2-52k-32)n \right)-\frac{1}{4}k^4-\frac{15}{4}k^3+12k^2+5k-\frac{7}{2} \end{eqnarray*}
\newpage
\section{Conclusion}
In this paper, the Zagreb indices and coindices of the line graph and line cut-vertex graph of the subdivision graph of the cycle-star graph are calculated.  For further research, it would be interesting to determine many other topological indices and coindices of cycle-star graphs using graph operators. For more details on graph operators, the readers are referred to \cite{18}. 
\makeatletter
\renewcommand{\@biblabel}[1]{[#1]\hfill}
\makeatother


\begin{thebibliography}{99}
\bibitem{1}{T. Do\v{s}li\'c, Vertex-weighted Weiner polynomials for composite graphs, \emph{Ars Math. Contemp}. \textbf{1} (2008) 66-80}.
\bibitem{2}{K. C. Das, I. Gutman, Some properties of the second Zagreb index, \emph{MATCH Commun. Math. Comput. Chem}. \textbf{52} (2004) 103-112}.
\bibitem{3}{B. Furtula, I. Gutman, M. Dehmer, On structure-sensitivity of degree-based topological indices, \emph{Appl. Math. Comput}. \textbf{219} (2013) 8973-8978}.
\bibitem{4}{I. Gutman, N. Trinajsti\'c N, Graph theory and molecular orbitals. Total $\pi$-electron energy of alternant hydrocarbons, \emph{Chem Phys Lett}. \textbf{17} (1972) 535-538}.
\bibitem{5}{I. Gutman, B. Ru\v{s}\v{c}i\'c, N. Trinajsti\'c, C. F. Wilcox, Graph theory and molecular orbitals. XII. Acyclic polyenes, \emph{J Chem Phys}. \textbf{62} (1975) 3399-3405}.
\bibitem{6}{I. Gutman, Degree-based topological indices, \emph{Croat. Chem, Acta}. \textbf{86} (2013) 351-361}.
\bibitem{7}{I. Gutman, K. C. Das, The first Zagreb index 30 years after, \emph{MATCH Commun. Math. Comput. Chem}. \textbf{50} (2004) 83-92}.
\bibitem{8}{I. Gutman, B. Furtula, Z. K. Vuki\'cevi\'c, G. Popivoda, On Zagreb Indices and Coindices, \emph{MATCH Commun. Math. Comput. Chem}. \textbf{74} (2015) 5-16}.
\bibitem{9}{I. Gutman, J.To\v{s}vi\'c, Testing the quality of molecular structure descriptors. Vertex-degree based topological indices, \emph{J. Serb. Chem. Soc}. \textbf{78} (2013) 805-810}.
\bibitem{10}{I. Gutman, On the origin of two degree-based topological indices, \emph{Bull. Acad. Serbe Sci. Arts (Cl. Sci. Math. Natur.)}. \textbf{146} (2014) 39-52}.
\bibitem{11}{F. Harary, \emph{Graph Theory}, Addison-Wesley, Reading, Mass (1969)}.
\bibitem{12}{V. R. Kulli, M. H. Muddebihal, On lict and litact graph of a graph, \emph{Proceeding of the Indian National Science Academy}. \textbf{41}(1975) 275-280}. 
\bibitem{13}{Monther R. Alfuraidan, Selvaraj Balachandran and Tomáš Vetr\'ik, General multiplicative Zagreb indices of unicyclic graphs, \emph{Carpathian Journal of Mathematics}. \textbf{37} (2021) 1-13}.
\bibitem{14}{S. Nikoli\'c, G. Kova\v{c}evi\'c, A. Mili\'cevi\'c, N. Trinajsti\'c, The Zagreb indices 30 years after, \emph{Croat Chem Acta}. \textbf{76} (2003) 113-124}.
\bibitem{15}{H. M. Nagesh, On the Zagreb indices of the line cut-vertex graph of subdivision of graphs, \emph{Bull. Int. Math. Virtual Inst}. \textbf{12} (2022) 17-26}.
\bibitem{16}{H. M. Nagesh, V. R. Girish, On the entire Zagreb indices of the line graph and line cut-vertex graph of the subdivision graph, \emph{Open Journal of Mathematical Sciences}. \textbf{4(2)} (2020) 470--475}
\bibitem{17}{H. M. Nagesh, V. R. Girish, B. Azghar Pasha, On the entire Zagreb indices of the line graph and line cut-vertex graph of the subdivision graph, \emph{Open Journal of Mathematical Sciences}. \textbf{12(1)} (2022) 159-167}
\bibitem{18}{E. Prisner, \emph{Graph dynamics}, Chapman \& Hall/CRC, \textbf{338} (1995)}.
\bibitem{19}{P. S. Ranjini, V. Lokesha, I. N. Kangul, On the Zabreb indices of the line graphs of the subdivision graphs, \emph{Appl. Math. Comput}. \textbf{210} (2011) 699-702}.
\bibitem{20}{J. Sedlar, Extremal unicyclic graphs with respect to additively weighted Harary index, \emph{Miskolic 
mathematical Notes}. \textbf{16(2)} (2013) 1-16}.
\bibitem{21}{Xuli Qi, Bo Zhou, Jiyong Li, Zagreb eccentricity indices of unicyclic graphs, \emph{Discrete Applied Mathematics}. \textbf{233(2)} (2017) 166-174}.
\bibitem{22}{Zikai Tang, Hechao Liu, Maximal hyper-Zagreb index of trees, unicyclic and bicyclic graphs with a given order and
matching number, \emph{Discrete Math. Lett}. \textbf{4} (2020) 11–18}.
\end{thebibliography}
\end{document}